\documentclass{birkart}
\usepackage{amssymb,mathrsfs,epsfig}

\renewcommand{\subjclassname}{\textup{2000} Mathematics Subject Classification}

 \newtheorem{thm}{Theorem}[section]
 \newtheorem{cor}[thm]{Corollary}
 \newtheorem{lem}[thm]{Lemma}
 \newtheorem{prop}[thm]{Proposition}
 \theoremstyle{definition}
 
 \theoremstyle{remark}
 
 \newtheorem*{ex}{Example}
 \numberwithin{equation}{section}

\def\csym{\,{}^c\!\sym}
\def\C{\mathbb C}
\def\Dom{\mathcal D}
\def\eps{\varepsilon}
\def\Gr{\mathrm{Gr}}
\def\K{\mathcal K}
\def\L{\mathscr L}
\def\open#1{\smash[t]{\overset{{}_{\circ}}{#1}{}}}
\def\R{\mathbb R}
\def\S{\mathscr S}
\def\Sing{\mathcal E}

\DeclareMathOperator{\bgres}{bg-res}

\DeclareMathOperator{\Diff}{Diff}
\DeclareMathOperator{\Ind}{ind}
\DeclareMathOperator{\LinSpan}{span}
\DeclareMathOperator{\res}{res}
\DeclareMathOperator{\spec}{spec}
\DeclareMathOperator{\sym}{ \sigma\!\!\!\sigma}

\begin{document}
%
\title[Rays of minimal growth]{On rays of minimal growth for \\
elliptic cone operators}
\author{Juan B. Gil}
\address{%
3000 Ivyside Park\\ Penn State Altoona\\
Altoona, PA 16601\\ U.S.A.}
\email{jgil@psu.edu}

\author{Thomas Krainer}
\address{%
Institut f\"ur Mathematik\\ Universit\"at Potsdam\\
D-14415 Potsdam\\ Germany}
\email{krainer@math.uni-potsdam.de}

\author{Gerardo A. Mendoza}
\address{%
Department of Mathematics\\ Temple University\\
Philadelphia, PA 19122\\ U.S.A.}
\email{gmendoza@math.temple.edu}

\subjclass{Primary 58J50; Secondary 35J70, 47A10}

\keywords{Resolvents, conical singularities, spectral theory}

\date{\today}

\begin{abstract}
We present an overview of some of our recent results on the existence of 
rays of minimal growth for elliptic cone operators and  two new results 
concerning the necessity of certain conditions for the 
existence of such rays.
\end{abstract}

\maketitle

\section{Introduction}

The aim of this article is twofold. On the one hand, we present an
overview of some of the results contained in \cite{GKM1,GKM2} on the
subject in the title, and of the geometric perspective we developed in the
course of the investigations leading to the aforementioned papers. We
illustrate the main ideas of our approach by means of examples concerning
Laplacians on a compact $2$-manifold. Already this simple situation
exhibits the structural richness and complexity of the general theory.

On the other hand, we offer some improvements, cf.
Theorems~\ref{GeometricNecAndSuff} and \ref{NecessityConstCoeff},
regarding necessary and sufficient conditions for a closed sector
$\Lambda\subset \C$ to be a sector of minimal growth for a certain class
of elliptic cone operators $A$ and for the associated model operator
$A_{\wedge}$.

Recall that a closed sector of the form 
\begin{equation}\label{Sector}
\Lambda = \{\lambda \in {\C}: \lambda = re^{i\theta} \; \text{for} \;
r \geq 0, \; \theta \in {\R}, \; |\theta - \theta_0| \leq a\}
\end{equation}
is called a \emph{sector of minimal growth} (or of maximal decay) for a 
closed operator
\[ A:\Dom\subset H\to H, \]
where $H$ is a Hilbert space and $\Dom$ is dense in $H$, 
if there is a constant $R>0$ such that $A-\lambda$ is invertible for 
every $\lambda\in \Lambda_R=\{\lambda\in\Lambda: |\lambda|\ge R\}$, and 
the resolvent $(A-\lambda)^{-1}$ satisfies either of the equivalent estimates  
\begin{equation}\label{NormEstimates}
 \bigl\|(A-\lambda)^{-1}\bigr\|_{\L(H)} \le C/|\lambda|, \quad
 \bigl\|(A-\lambda)^{-1}\bigr\|_{\L(H,\Dom)} \le C 
\end{equation}
for some $C>0$ and all $\lambda\in\Lambda_R$.

We are interested in cone operators on smooth manifolds with boundary.
Specifically, let $M$ be a smooth $n$-manifold with boundary $Y=\partial
M$ and let $E\to M$ be a Hermitian vector bundle over $M$. Fix a defining
function $x$ for $Y$. A differential \emph{cone operator} of order $m$
acting on sections of $E\to M$ is an operator of the form $A= x^{-m}P$
with $P$ in the class $\Diff_b^m(M;E)$ of totally characteristic
differential operators of order $m$, cf. Melrose \cite{RBM1}. We write
$A\in x^{-m}\Diff_b^m(M;E)$. More explicitly, in the interior $\open M$
of $M$, $A$ is a differential operator with smooth coefficients, and near
the boundary, in local coordinates $(x,y)\in (0,\eps)\times Y$, it is of
the form 
\begin{equation}\label{cone-operator} 
  A=x^{-m}\sum_{k+|\alpha|\le m} a_{k\alpha}(x,y)(xD_x)^k D_y^\alpha 
\end{equation} 
with coefficients $a_{k\alpha}$ smooth up to $x=0$; here 
$D_x=-i\partial/\partial x$ and likewise $D_{y_j}$. We will say that $A$
(or $P$) has \emph{coefficients independent of $x$ near $Y$}, if the
coefficients $a_{k\alpha}$ in \eqref{cone-operator} do not depend on $x$
(this notion depends on the choice of tubular neighborhood map, defining
function $x$, and connection on $E$. For a precise definition see \cite{GKM1}).

This paper consists of 5 sections. In Section 2 we review some basic 
properties of cone operators while in Section 3 we discuss the associated
model operators.  The new results on rays of minimal growth can be found in 
Sections 4 and 5.  Apart from the works explicitly cited in the text, 
our list of references contains additional items referring to related 
works on resolvents and rays of minimal growth for elliptic operators.

\section{Preliminaries on cone operators}
\label{sec-Preliminaries}

Let $A$ be a differential cone operator. As introduced in \cite{GKM1}, the 
principal symbol of $A$, $\csym(A)$, is defined on the $c$-cotangent 
${}^cT^*M$ of $M$ rather than on the cotangent itself. Over $\open M$ it is 
essentially the usual principal symbol, and equal to
\begin{equation*}
 \sum_{k+|\alpha|= m} a_{k\alpha}(x,y)\xi^k \eta^\alpha
\end{equation*}
near the boundary $Y$, see \eqref{cone-operator}.  

\begin{ex}
Let $M$ be a compact $2$-manifold with boundary $Y=S^1$. Let $g_{Y}(x)$ be
a smooth family of Riemannian metrics on $S^1$ such that $g_{Y}(0)$ is the 
standard metric, $dy^2$.  We equip $M$ with a ``cone metric'' $g$ 
that near $Y$ takes the form $g=dx^2 + x^2 g_{Y}(x)$ ($g$ is a regular 
Riemannian metric in the interior of $M$).  Then, near $Y$, the 
Laplace-Beltrami operator $\Delta$ has the form 
\begin{equation}\label{Laplacian}
 x^{-2}\big((xD_x)^2 + a(x,y)(xD_x) + \Delta_{Y}(x)\big),  
\end{equation}
where $a(x,y)$ is a smooth function with $a(0,y)=0$ and $\Delta_Y(x)$  
is the nonnegative Laplacian on $S^1$ associated with $g_{Y}(x)$. 
In this case, near the boundary, we have
\[ \csym(\Delta) = \xi^2 + \sym(\Delta_Y(x)). \]
\end{ex}

\subsection*{Ellipticity and boundary spectrum}

An operator $A\in x^{-m}\Diff_b^m(M;E)$ is \emph{$c$-elliptic} if 
$\csym(A)$ is invertible on ${}^cT^*M\backslash 0$. Moreover, the  
family $A-\lambda$ is said to be \emph{$c$-elliptic with parameter}  
$\lambda\in\Lambda\subset\C$  if $\csym(A)-\lambda$ is invertible on
$({}^cT^*M\times\Lambda) \backslash 0$.

Associated with $A=x^{-m}P$ there is an operator-valued polynomial
\begin{equation*}
 \C\ni \sigma\mapsto \hat P(\sigma)\in \Diff^m(Y;E|_Y)
\end{equation*}
called the \emph{conormal symbol} of $P$ (and of $A$). 
If we write $A$ as in \eqref{cone-operator}, then
\begin{equation*}
 \hat P(\sigma) = \sum_{k+|\alpha|\le m}
 a_{k\alpha}(0,y)\sigma^k D_y^\alpha.
\end{equation*}
If $A$ is $c$-elliptic, then $\hat P(\sigma)$ is invertible for all
$\sigma\in\C$ except a discrete set $\spec_b(A)$, the \emph{boundary 
spectrum of $A$}, cf. \cite{RBM1};  $\hat P(\sigma)$ is a holomorphic 
family of elliptic operators on $Y$ and $\sigma\to \hat P(\sigma)^{-1}$ is 
a meromorphic operator-valued function on $\C$.

\begin{ex}
The Laplacian \eqref{Laplacian} is clearly $c$-elliptic. If $y$ is the 
angular variable on $S^1$, then 
\begin{equation*}
 \hat P(\sigma) = \sigma^2 + \Delta_Y(0) =\sigma^2+D_y^2,
\end{equation*}
and the boundary spectrum of $\Delta$ is given by
\begin{equation*}
 \spec_b(\Delta) = \{\pm ik: k\in\mathbb{N}_0\}.
\end{equation*}
\end{ex}

\subsection*{Closed extensions}

Let $\mathfrak{m}$ be a positive $b$-density on $M$, that is, 
$x\mathfrak{m}$ is a smooth everywhere positive density on $M$.  
Let $L^2_b(M;E)$ be the $L^2$ space of sections of $E$ with respect to the
Hermitian form on $E$ and the density $\mathfrak{m}$.  Consider $A$ 
initially defined on $C_0^\infty(\open M;E)$ and look at it as an unbounded 
operator on the Hilbert space 
\begin{equation*}
x^{-m/2}L^2_b(M;E) = L^2(M;E;x^{2m} \mathfrak m).
\end{equation*}
The particular weight $x^{-m/2}$ is just a convenient normalization and 
represents no loss.  If we are interested in $A$ on $x^{\mu}L^2_b(M;E)$ for 
$\mu\in\R$, we can base all our analysis on the space $x^{-m/2}L^2_b(M;E)$ 
by considering the operator $x^{-\mu-m/2} A\,x^{\mu+m/2}$.

Typically, $A$ has a large class of closed extensions
\begin{equation}\label{Ext}
A_{\Dom}:\Dom\subset x^{-m/2}L^2_b(M;E)\to x^{-m/2}L^2_b(M;E).
\end{equation}
There are two canonical closed extensions, namely the ones with domains
\begin{align*}
\Dom_{\min}(A) &= \text{closure of } C^\infty_0(\open M;E) 
  \text{ with respect to } \|\cdot\|_A, \\
\Dom_{\max}(A) &=\{u\in x^{-m/2}L^2_b(M;E): Au\in x^{-m/2}L^2_b(M;E)\},
\end{align*}
where $\|u\|_A=\|u\|+\|Au\|$ is the graph norm in $\Dom_{\max}(A)$.
Both domains are dense in $x^{-m/2}L^2_b(M;E)$, and for any closed extension 
\eqref{Ext},
\[ \Dom_{\min}(A)\subseteq \Dom\subseteq \Dom_{\max}(A). \]
Let
\begin{equation*}
  \mathfrak D(A)=\{\Dom\subset\Dom_{\max}(A): \Dom \text{ is a vector space
and } \Dom_{\min}(A)\subset\Dom\}.
\end{equation*}
The elements of $\mathfrak D(A)$ are in one-to-one correspondence with the
subspaces of $\Dom_{\max}(A)/\Dom_{\min}(A)$.  If the operator $A$ is fixed 
and there is no possible ambiguity, we will omit $A$ from the notation and 
will write simply $\Dom_{\min}$, $\Dom_{\max}$, and $\mathfrak D$.

\begin{thm}[Lesch \cite{Le97}] \label{LeschTheorem}
If $A\in x^{-m}\Diff_b^m(M;E)$ is $c$-elliptic, then
\[ \dim \Dom_{\max}/\Dom_{\min}<\infty \]
and all closed extensions of $A$ are Fredholm.
Moreover, 
\begin{equation}\label{RelIndexA}
\Ind A_{\Dom}= \Ind A_{\Dom_{\min}}+ \dim\Dom/\Dom_{\min}. 
\end{equation}
\end{thm}

Modulo $\Dom_{\min}$, the elements of $\Dom_{\max}$ are determined by their 
asymptotic behavior near the boundary of $M$. The structure of these 
asymptotics depends on the conormal symbols of $A$ and on the part of 
$\spec_b(A)$ in the strip $\{|\Im\sigma|<m/2\}$. More details will be
discussed in the next section. 

\begin{cor}
If $A$ is $c$-elliptic and symmetric (formally selfadjoint), then
\begin{equation*}
 \Ind A_{\Dom_{\max}} = -\Ind A_{\Dom_{\min}} \;\;\text{and}\;\;
 \Ind A_{\Dom_{\min}} = -\frac12 \dim\Dom_{\max}/\Dom_{\min}.
\end{equation*}
\end{cor}

\begin{ex}
Consider the cone Laplacian $\Delta$, cf. \eqref{Laplacian}.
Then \eqref{EqualDimension} and \eqref{QuotientLaplacian} imply
\begin{equation*}
 \dim\Dom_{\max}(\Delta)/\Dom_{\min}(\Delta)=2  
\end{equation*}
and thus, by the previous corollary,
\begin{equation}\label{IndexLaplacian}
 \Ind \Delta_{\min}=-1 \;\;\text{ and }\;\; \Ind \Delta_{\max}=1. 
\end{equation}
\end{ex}

If $A\in x^{-m}\Diff_b^m(M;E)$ is $c$-elliptic, the embedding
$\Dom_{\max}\hookrightarrow x^{-m/2}L^2_b(M;E)$ is compact. Therefore,
for every $\Dom\in\mathfrak D$ and $\lambda\in\C$, the operator 
$A_{\Dom}-\lambda$ is also Fredholm with $\Ind (A_{\Dom}-\lambda)=
\Ind A_{\Dom}$. Consequently, if $\spec(A_\Dom)\not=\C$, then we 
necessarily have $\Ind A_\Dom=0$. For this reason, we will primarily be
interested in the set of domains
\begin{equation}\label{GrasmannianA}
 \mathfrak G=\{\Dom\in \mathfrak D: \Ind A_\Dom = 0\}
\end{equation}
which is empty unless $\Ind A_{\Dom_{\min}}\le 0$ and $\Ind
A_{\Dom_{\max}}\ge 0$. Let $d''=-\Ind A_{\Dom_{\min}}$. 

Using that the map $\mathfrak D\ni \Dom\mapsto \Dom/\Dom_{\min}$ is a 
bijection, we identify $\mathfrak G$ with the complex Grassmannian of 
$d''$-dimensional subspaces of $\Dom_{\max}/\Dom_{\min}$.

\begin{ex} 
For $\Delta$ we have
\begin{equation*}
 \mathfrak G(\Delta) \cong \mathbb{CP}^1 = S^2.
\end{equation*}
Note that by \eqref{RelIndexA} and \eqref{IndexLaplacian}, $\Ind\Delta_{\Dom}=0$ if and only if $\dim\Dom/\Dom_{\min}=1$. 
\end{ex}

We finish this section with the following proposition that gives a first 
glimpse of the complexity of the spectrum of elliptic cone operators.
\begin{prop}
If $A$ is $c$-elliptic and $\dim \mathfrak G>0$, then for any $\lambda\in\C$ 
there is a domain $\Dom\in\mathfrak G$ such that $\lambda\in\spec(A_\Dom)$.
If, in addition, $A$ is symmetric on $\Dom_{\min}$, then for any 
$\lambda\in\R$ there is a $\Dom\in\mathfrak G$ such that $A_\Dom$ is 
selfadjoint and $\lambda\in\spec(A_\Dom)$.
\end{prop}
A proof is given in \cite[Propositions~5.7 and ~6.7]{GKM1}.  
A surprising consequence of the second statement is that for any arbitrary 
negative number $\lambda$ there is always a selfadjoint extension of $A$ 
having $\lambda$ as eigenvalue, even if $A$ is positive on $\Dom_{\min}$.

\section{The model operator} \label{sec-ModelOperator}

Let $A\in x^{-m}\Diff_b^m(M;E)$ be $c$-elliptic. The \emph{model operator 
$A_\wedge$} associated with $A$ is an operator on $N_+Y$, the closed inward 
normal bundle of $Y$, that in local coordinates takes the form  
\begin{equation*}
 A_\wedge=x^{-m}\sum_{k+|\alpha|\le m}
 a_{k\alpha}(0,y)(xD_x)^k D_y^\alpha,
\end{equation*}
if $A$ is written as in \eqref{cone-operator}. A Taylor expansion in $x$ (at $x=0$) of the coefficients of the operator $A$ induces a decomposition
\begin{equation}\label{bOperatorTaylor}
 x^m A = \sum_{k=0}^{N-1} P_k x^k + x^N \tilde P_N
 \quad\text{for every } N\in\mathbb{N}, 
\end{equation}
where each $P_k$ has coefficients independent of $x$ near $Y$.  Thus the 
model operator can be written, near $Y$, as $A_\wedge =x^{-m}P_0$.  In other 
words, $A_\wedge$ can be thought of as the ``most singular'' part of $A$.

We trivialize $N_+Y$ as $Y^\wedge=[0,\infty)\times Y$. 
The operator $A_\wedge\in x^{-m}\Diff_b^m(Y^\wedge;E)$ acts 
on $C^\infty_0(\open Y^\wedge;E)$ and can be extended as a densely defined 
closed operator in $x^{-m/2}L^2_b(Y^\wedge;E)$. 
The space $L^2_b(Y^\wedge;E)$ is the $L^2$ space with respect to 
a density of the form $\frac{dx}{x}\otimes \pi^*\mathfrak{m}_Y$ and 
the canonically induced Hermitian form on $\pi^*(E|_Y)$, where 
$\pi : Y^{\wedge} \to Y$ is the projection on the factor $Y$. The density
$\mathfrak{m}_Y$ is related to $\mathfrak{m}$ and, by abuse of notation, 
we denote $\pi^*(E|_Y)$ by $E$, cf. \cite{GKM1}.
Again, there are two canonical domains $\Dom_{\wedge,\min}$ and 
$\Dom_{\wedge,\max}$ and we denote by $\mathfrak D_\wedge$ the set of 
subspaces of $\Dom_{\wedge,\max}$ that contain $\Dom_{\wedge,\min}$. 
There is a natural (and useful) linear isomorphism
\begin{equation*}
 \theta: \Dom_{\max}/\Dom_{\min} \to \Dom_{\wedge,\max}/\Dom_{\wedge,\min},
\end{equation*}
cf. Section~\ref{sec-Necessity}.  As a consequence we have 
\begin{equation}\label{EqualDimension}
 \dim \Dom_{\wedge,\max}/\Dom_{\wedge,\min} = \dim \Dom_{\max}/\Dom_{\min}
\end{equation}
which by Theorem~\ref{LeschTheorem} is finite. It is known (cf. 
Lesch~\cite{Le97}) that $\Dom_{\wedge,\max}/\Dom_{\wedge,\min}$ is isomorphic 
to a finite dimensional space $\Sing_{\wedge,\max}$ consisting of functions 
of the form
\begin{equation}\label{SingularFunctions} 
  \varphi=\sum_{\substack{\sigma\in\spec_b(A) \\|\Im\sigma|<m/2}} 
  \left(\sum_{k=0}^{m_\sigma} c_{\sigma,k}(y) \log^k x\right) x^{i\sigma} 
\end{equation} 
where $c_{\sigma,k}\in C^\infty(Y;E)$. More precisely, for every 
$u\in\Dom_{\wedge,\max}$ there is a function $\varphi\in \Sing_{\wedge,\max}$ 
such that $u(x,y)-\omega(x)\varphi(x,y)\in \Dom_{\wedge,\min}$ for some 
(hence any) cut-off function $\omega \in C_0^\infty([0,1))$, $\omega=1$ near 
$0$. The function $\varphi$ is uniquely determined by the equivalence class 
$u+\Dom_{\wedge,\min}$. 

We identify
\[ \Sing_{\wedge,\max} = \Dom_{\wedge,\max}/\Dom_{\wedge,\min} \]
and let
\begin{equation*}
  \pi_{\wedge,\max}:\Dom_{\wedge,\max}\to \Sing_{\wedge,\max}
\end{equation*}
be the canonical projection.

Contrary to the situation in Theorem~\ref{LeschTheorem}, the closed 
extensions of $A_\wedge$ do not need to be Fredholm. However, 
if $A-\lambda$ is $c$-elliptic with parameter, then the canonical 
extensions $A_{\wedge,\min}-\lambda$ and $A_{\wedge,\max}-\lambda$ are both 
Fredholm for $\lambda\not=0$, cf. \cite[Remark~5.26]{GKM2}.
Moreover, we have
\begin{equation}\label{IndexAwedgemin}
 \Ind(A_{\wedge,\min}-\lambda)=\Ind A_{\Dom_{\min}},
\end{equation}
cf. Corollary~{5.35} in \cite{GKM2}.

\begin{ex}
On $Y^\wedge=[0,\infty)\times S^1$ with the cone metric $dx^2+x^2 dy^2$, 
the Laplace-Beltrami operator is given by
\begin{equation}\label{ModelLaplacian}
 \Delta_\wedge = x^{-2}\big((xD_x)^2 + \Delta_Y\big),
\end{equation}
where $\Delta_Y$ is the nonnegative Laplacian on $S^1$. $\Delta_\wedge$ is 
precisely the model operator associated with the cone Laplacian $\Delta$ 
discussed in the previous section, cf. \eqref{Laplacian}. 
It is easy to check that for any cut-off function $\omega\in
C_0^\infty([0,1))$, the functions 
\[ \omega(x)\cdot 1,\; \Delta_\wedge(\omega(x)\cdot 1), \;
   \omega(x)\log x,\, \text{ and } \Delta_\wedge(\omega(x)\log x) \]
are all in the space $x^{-1}L^2_b(Y^\wedge)$. Thus 
$\omega(x)\cdot 1$ and $\omega(x)\log x$ are elements of $\Dom_{\wedge,\max}$.
In fact,
\begin{equation}\label{QuotientLaplacian}
 \Sing_{\wedge,\max}=\LinSpan\{1,\log x\}.
\end{equation}
Observe that $\Delta-\lambda$ is $c$-elliptic with parameter
$\lambda\in \C\backslash \R_+$ and therefore the closed extensions of 
$\Delta_\wedge-\lambda$ are Fredholm for every $\lambda\in \C\backslash \overline{\R}_+$.
\end{ex}

The model operator has a dilation/scaling property that can be 
exploited to analyze its closed extensions and their resolvents from a 
geometric point of view. In order to describe this property we first 
introduce the one-parameter group of isometries
\begin{equation*}
\R_+\ni \varrho \mapsto
\kappa_\varrho:x^{-m/2}L^2_b(Y^\wedge;E)\to x^{-m/2}L^2_b(Y^\wedge;E)
\end{equation*}
which on functions is defined by
\begin{equation}\label{kapparho}
  (\kappa_\varrho f)(x,y)= \varrho^{m/2} f(\varrho x,y).
\end{equation}
It is easily verified that the operator $A_\wedge$ satisfies the relation
\begin{equation*}
  \kappa_\varrho A_\wedge =\varrho^{-m} A_\wedge\kappa_{\varrho}.
\end{equation*}
This implies 
\begin{equation}\label{kappaHomogeneous}
  A_\wedge-\lambda = \varrho^{m} \kappa_\varrho (A_\wedge-\lambda/\varrho^m) 
  \kappa_{\varrho}^{-1} 
\end{equation}
for every $\varrho>0$ and $\lambda\in\C$. This homogeneity property,
called \emph{$\kappa$-homogeneity}, will be used systematically to describe 
the closed extensions of $A_\wedge$ with nonempty resolvent sets.

It is convenient to introduce the set
\begin{equation*}
 \bgres A_\wedge=\{\lambda\in\C: A_{\wedge,\min}-\lambda \text{ is
 injective and }  A_{\wedge,\max}-\lambda \text{ is surjective}\},
\end{equation*}
the \emph{background resolvent set} of $A_\wedge$, cf.~\cite{GKM1}.
\begin{lem}[Lemma 7.3 in \cite{GKM1}]
If $\lambda\in \bgres A_\wedge$ and $\Dom\in \mathfrak D_\wedge$, then 
$A_{\wedge,\Dom} -\lambda$ is Fredholm. The set $\bgres A_\wedge$ is a 
disjoint union of open sectors, 
\begin{equation*}
 \bgres A_\wedge = \bigcup_{\alpha\in \mathfrak I\subset \mathbb N} 
 \open \Lambda_\alpha.
\end{equation*}
\end{lem}
This lemma follows immediately from \eqref{kappaHomogeneous}.

For $\lambda\in \bgres A_\wedge$ and $\Dom\in \mathfrak D_\wedge$ we have
\begin{equation}\label{RelIndexAwedge}
 \Ind(A_{\wedge,\Dom}-\lambda)=
 \Ind(A_{\wedge,\min}-\lambda)+\dim\Dom/\Dom_{\wedge,\min}.
\end{equation}
Moreover, the map
\begin{equation*}
 \open\Lambda_\alpha\ni \lambda\mapsto \Ind(A_{\wedge,\Dom}-\lambda)
\end{equation*}
is constant since the embedding $\Dom\hookrightarrow 
x^{-m/2}L^2_b(Y^\wedge;E)$ is continuous.
Now, in analogy with \eqref{GrasmannianA} we define 
\begin{equation*}
 \mathfrak G_{\wedge,\alpha}=\{\Dom\in\mathfrak D_\wedge:
 \Ind(A_{\wedge,\Dom}-\lambda)=0 \text{ for } \lambda\in
 \open\Lambda_\alpha\}
\end{equation*}
and let $d''_\alpha= -\Ind(A_{\wedge,\min}-\lambda)$ for $\lambda\in 
\open\Lambda_\alpha$.  We identify $\mathfrak G_{\wedge, \alpha}$ with the 
complex Grassmannian of $d''_\alpha$-dimensional subspaces of 
$\Sing_{\wedge,\max}$. 

The canonical domains $\Dom_{\wedge,\min}$ and $\Dom_{\wedge,\max}$ are 
both $\kappa$-invariant. Thus the group action $\kappa_{\varrho}$ induces
an action on $\Sing_{\wedge,\max}$. In general, $\kappa_{\varrho}$ does 
not preserve the elements of $\mathfrak D_\wedge$. In fact, the set of 
$\kappa$-invariant domains in $\mathfrak D_\wedge$ is an analytic variety 
because it consists of the stationary points of a holomorphic flow,
cf. Section~{7} in \cite{GKM1}.  
To better analyze the resolvents of the closed extensions of $A_\wedge$
over the open sector $\open \Lambda_\alpha$, we will consider the manifold
$\mathfrak G_{\wedge,\alpha}$ together with the flow generated by the
induced action of $\kappa_{\varrho}$ given by 
$\kappa_{\varrho}(\Dom/\Dom_{\wedge,\min})=
 \kappa_{\varrho}(\Dom)/\Dom_{\wedge,\min}$.

\begin{ex}
The background resolvent set of $\Delta_\wedge$ is the open sector
$\C\backslash \overline{\R}_+$; this is easily seen after noting that 
$\Delta_\wedge$ is the standard Laplacian in $\R^2$ written in polar 
coordinates. Moreover, since $\Ind(\Delta_{\wedge,\min}-\lambda)= -1$ 
for every $\lambda\in \C\backslash \overline{\R}_+$, we have that 
$\Dom\in\mathfrak D_{\wedge}$ belongs to $\mathfrak G_\wedge$ if and 
only if $\dim \Dom/\Dom_{\wedge,\min}=1$. Thus
\begin{equation*}
 \mathfrak G_\wedge \cong \mathbb{CP}^1 = S^2.
\end{equation*}
We identify $\Sing_{\wedge,\max}$ with $\Dom_{\wedge,\max}/\Dom_{\wedge,\min}$ 
and use \eqref{QuotientLaplacian} to write
\begin{equation*}
 \Sing_{\wedge,\max} = \LinSpan\left\{1,\log x\right\}.
\end{equation*}
For $\Dom\in \mathfrak G_\wedge$ we then have
\begin{equation}\label{FiniteDomainLaplace}
 \pi_{\wedge,\max}\Dom = \LinSpan\left\{\zeta_0\cdot 1+\zeta_1\log x\right\} 
 \;\text{ for some } \zeta_0,\zeta_1\in\C, \, (\zeta_0,\zeta_1)\not=0.
\end{equation} 
Hence, with $\kappa$ as defined in \eqref{kapparho}, we get
\begin{equation}\label{FinitekappaD}
 \pi_{\wedge,\max}\kappa_{\!\varrho}^{-1}\Dom 
 =\LinSpan\{(\zeta_0-\zeta_1\log\varrho)\cdot 1 + \zeta_1\log x\}.
\end{equation}
Clearly, the only $\kappa$-invariant domain in $\mathfrak G_\wedge$ is the 
domain $\Dom_{F}$ such that
\[ \pi_{\wedge,\max}\Dom_{F}=\LinSpan\{1\}; \]
$\Dom_F$ is precisely the domain of the Friedrichs extension of 
$\Delta_\wedge$, cf. \cite{GiMe01}. Every domain $\Dom\in
\mathfrak G_\wedge$ with $\zeta_1\not=0$ in \eqref{FiniteDomainLaplace}
generates a nontrivial orbit as given by \eqref{FinitekappaD}.
In order to describe the flow of $\kappa$ on these nonstationary points,
rewrite \eqref{FiniteDomainLaplace} as
\begin{equation*}
 \pi_{\wedge,\max}\Dom = 
 \LinSpan\left\{\tfrac{\zeta_0}{\zeta_1}\cdot 1+\log x\right\}. 
\end{equation*}
Then the projection to $\Sing_{\wedge,\max}$ of the dilation 
$\kappa_{\varrho}^{-1}\Dom$ is given by
\begin{equation}\label{FinitekappaD2}
 \pi_{\wedge,\max}\kappa_{\!\varrho}^{-1}\Dom 
 =\LinSpan\left\{\big(\tfrac{\zeta_0}{\zeta_1}-\log\varrho\big)\cdot 1 
  + \log x\right\}.
\end{equation}
If $[\zeta_0:\zeta_1]\in\mathbb{CP}^1$ is the point corresponding to $\Dom$, 
then $\kappa_{\varrho}^{-1}\Dom$ is represented by
$[\zeta_0-\zeta_1\log\varrho:\zeta_1]$. In other words, in the situation 
at hand, the flow generated by $\kappa$ on $\mathfrak G_\wedge \cong 
\mathbb{CP}^1$ consists of curves that in projective coordinates are lines 
parallel to the real axis, see Figure \ref{figOrbits}. 
\begin{figure}[!ht]
\begin{picture}(100,130)
\put(-30,-10){\epsfig{file=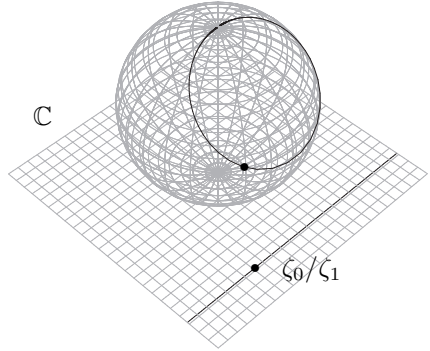,width=2.2in}}
\put(64,21){\circle*{3}} \put(74,18){$\zeta_0/\zeta_1$}
\put(60,59){\circle*{3}} \put(-20,75){$\C$}
\end{picture}
\caption{Orbit in $\mathfrak{G}_\wedge(\Delta_\wedge)$ generated 
by $\Dom\leftrightarrow \zeta_0/\zeta_1\in\C$.} \label{figOrbits}
\end{figure}

Observe that the Friedrichs extension corresponds to the point $[1:0]\in
\mathbb{CP}^1$. Using 
\begin{equation*}
 \pi_{\wedge,\max}\kappa_{\!\varrho}^{-1}\Dom 
 =\LinSpan\left\{1+\tfrac{\zeta_1}{\zeta_0-\zeta_1\log\varrho}
  \log x\right\}\;\text{ if } \varrho\not=e^{\zeta_0/\zeta_1}
\end{equation*}
we see that
\begin{equation}\label{DomainConvergence}
 \pi_{\wedge,\max}\kappa_{\!\varrho}^{-1}\Dom \to 
 \LinSpan\{1\} = \pi_{\wedge,\max}\Dom_F
 \;\text{ as $\varrho\to\infty$ or $\varrho\to 0$.} 
\end{equation}

Let $\Dom_0\in\mathfrak G_\wedge$ be such that $\pi_{\wedge,\max}\Dom_0 
=\LinSpan\{\log x\}$. This domain gives a selfadjoint extension of 
$\Delta_\wedge$ which on the sphere corresponds to the point $[0:1]$. The 
circle consisting of the orbit of $\Dom_0$ together with $\Dom_F$ is the 
set of domains of selfadjoint extensions of $\Delta_\wedge$.
\end{ex}

\section{Ray conditions on the model cone}

In this section we will discuss the existence of sectors of minimal growth
for the model operator $A_\wedge\in x^{-m}\Diff^m_b(Y^\wedge;E)$ 
associated with a $c$-elliptic cone operator. We fix a component
$\open\Lambda_\alpha$ of $\bgres A_\wedge$ and let 
$\mathfrak G_{\wedge}=\mathfrak G_{\wedge,\alpha}$. 

Let $\Lambda$ be a closed sector such that $\Lambda\backslash 0\subset 
\open\Lambda_\alpha$, cf. \eqref{Sector}, let $\res A_{\wedge,\Dom}$
be the resolvent set of $A_{\wedge,\Dom}$.
The $\kappa$-invariant domains are the simplest domains to analyze. 

\begin{prop}[Proposition~{8.4} in \cite{GKM1}]
Suppose $\Dom\in\mathfrak G_{\wedge}$ is $\kappa$-invariant. If there 
exists $\lambda_0\in\open\Lambda_\alpha$ such that $A_{\wedge,\Dom}-
\lambda_0$ is invertible, then $\open\Lambda_\alpha\subset
\res A_{\wedge,\Dom}$ and $\Lambda$ is a sector of minimal growth for 
$A_{\wedge,\Dom}$.
\end{prop}

If $\Dom\in\mathfrak G_{\wedge}$ is not $\kappa$-invariant, the situation
is more complicated. Nonetheless, in \cite{GKM1} we found a condition 
necessary and sufficient for a sector $\Lambda$ to be a sector of minimal 
growth for $A_{\wedge,\Dom}$. This condition is expressed in terms of 
finite dimensional spaces and projections that we proceed to discuss briefly. 

For $\lambda\in\bgres A_\wedge$ we let
\begin{equation*}
 \K_{\wedge,\lambda}=\ker(A_{\wedge,\max}-\lambda).
\end{equation*} 
Then 
\begin{equation*}
 \res A_{\wedge,\Dom}=\bgres A_\wedge\cap \{\lambda:
 \K_{\wedge,\lambda}\cap \Dom=0\}, 
\end{equation*}
and for $\lambda\in \res A_{\wedge,\Dom}$ we have
\begin{equation}\label{Dmax=KplusD}
 \Dom_{\wedge,\max} = \K_{\wedge,\lambda}\oplus \Dom.
\end{equation}
Projecting on $\Sing_{\wedge,\max}$, this direct sum induces the decomposition
\begin{equation}\label{EwedgeDecomposition}
 \Sing_{\wedge,\max}= \pi_{\wedge,\max}\K_{\wedge,\lambda}
 \oplus \pi_{\wedge,\max}\Dom,
\end{equation}
and the projection on $\pi_{\wedge,\max} \K_{\wedge,\lambda}$
according to \eqref{EwedgeDecomposition} is given by the map
\begin{equation}\label{KmaxProjection}
\begin{gathered}
 \hat\pi_{\K_{\wedge,\lambda}, \Dom}:
 \Sing_{\wedge,\max} \to\Sing_{\wedge,\max} \\
 (u+\Dom_{\wedge,\min}) \mapsto \pi_{\K_{\wedge,\lambda}, \Dom}u 
 + \Dom_{\wedge,\min},
\end{gathered}
\end{equation}
where $\pi_{\K_{\wedge,\lambda},\Dom}$ is the projection on 
$\K_{\wedge,\lambda}$ according to \eqref{Dmax=KplusD}.

\bigskip
The following theorem gives a condition on the operator norm of
\eqref{KmaxProjection} for a sector $\Lambda$ to be a sector of minimal
growth for $A_{\wedge,\Dom}$. Define
\begin{equation*}
 \|u\|^2_\lambda=\|u\|^2+|\lambda|^{-2}\|A_\wedge u\|^2
\end{equation*}
for $\lambda\not=0$ and $u\in\Dom_{\wedge,\max}$.

\begin{thm}\label{NecAndSuffWedge}
Let $\Dom\in \mathfrak G_{\wedge}$, let $\Lambda$ be a closed sector 
with $\Lambda\backslash 0\subset \open\Lambda_\alpha$. Then $\Lambda$ is
a sector of minimal growth for $A_{\wedge,\Dom}$ if and only if there are
$C$, $R>0$ such that $\Lambda_R\subset \res A_{\wedge,\Dom}$, and
\begin{equation}\label{NewFormWedgeCond}
 \big\|\hat\pi_{\K_{\wedge,\lambda},\Dom}
 \big\|_{\L(\Sing_{\wedge,\max},\|\cdot\|_\lambda)}\leq C 
 \;\text{ for } \lambda\in \Lambda_R,
\end{equation}
where $\hat\pi_{\K_{\wedge,\lambda},\Dom}$ is the projection 
\eqref{KmaxProjection}.
\end{thm}

This theorem is a rephrasing of \cite[Theorem~{8.7}]{GKM1}. There 
the condition \eqref{NewFormWedgeCond} appears in the equivalent form  
\begin{equation}\label{NecAndSuffWedgeCond} \tag{\ref{NewFormWedgeCond}$'$}
  \big\|\hat\pi_{\K_{\wedge,\hat\lambda},
  \kappa_{|\lambda|^{1/m}}^{-1}\Dom}
  \big\|_{\L(\Sing_{\wedge,\max})}\leq C \;\text{ for } \lambda\in \Lambda_R,
\end{equation}
where $\hat\lambda=\lambda/|\lambda|$, and $\hat\pi_{\K_{\wedge,\hat\lambda}, 
\kappa_{|\lambda|^{1/m}}^{-1}\Dom}$ is the projection on 
$\K_{\wedge,\hat\lambda}$ induced (following the steps 
\eqref{Dmax=KplusD}--\eqref{KmaxProjection}) by the direct sum
\begin{equation}\label{Dmax=kappaKplusD}
 \Dom_{\wedge,\max}= 
 \K_{\wedge,\varrho^{-m}\lambda}\oplus \kappa_\varrho^{-1}\Dom 
\end{equation}
for $\lambda\in\res A_{\wedge,\Dom}$ and $\varrho>0$. This decomposition
is a consequence of \eqref{Dmax=KplusD} and the $\kappa$-invariance of 
$\Dom_{\wedge,\max}$, as follows.  First, the $\kappa$-homogeneity of 
$A_\wedge-\lambda$, cf. \eqref{kappaHomogeneous}, implies 
\begin{equation*}
 \kappa_\varrho^{-1}(\K_{\wedge,\lambda})= 
 \K_{\wedge,\varrho^{-m}\lambda} \;\text{ for } \varrho>0.
\end{equation*}
Furthermore, if $\Dom\in\mathfrak G_\wedge$ and $\lambda\in\bgres
A_\wedge$, then 
\begin{equation*}
 \varrho^{-m}\lambda\in \res A_{\wedge,\kappa_\varrho^{-1}\Dom} 
 \Longleftrightarrow 
 \lambda\in \res A_{\wedge,\Dom}. 
\end{equation*}
In particular, 
\[ \K_{\wedge,\varrho^{-m}\lambda}\cap \kappa_\varrho^{-1}\Dom=\{0\}
   \Longleftrightarrow \K_{\wedge,\lambda}\cap \Dom=\{0\}, \]
as claimed.

The equivalence of \eqref{NewFormWedgeCond} and \eqref{NecAndSuffWedgeCond} 
follows immediately from the identity 
\begin{equation*}
  \kappa_{|\lambda|^{1/m}}^{-1} \pi_{\K_{\wedge,
  \lambda},\Dom}\, \kappa_{|\lambda|^{1/m}}= \pi_{\K_{\wedge,\hat\lambda},
  \kappa_{|\lambda|^{1/m}}^{-1}\Dom}
\end{equation*}
using the relation \eqref{kappaHomogeneous} and the fact that $\kappa$ is 
an isometry on $x^{-m/2}L^2_b$. The virtue of \eqref{NecAndSuffWedgeCond} 
is that the norm is fixed, while the advantage of \eqref{NewFormWedgeCond}
lies in that it gives a more explicit 
dependence on $\lambda$ and deals with a projection on a subspace of 
$\Sing_{\wedge,\max}$ with {\em fixed} complement $\pi_{\wedge,\max}\Dom$.

In \cite[Corollary {8.22}]{GKM1} it is proved that $\Lambda$ is a 
sector of minimal growth for $A_{\wedge,\Dom}$ if and only if there are 
constants $C$, $R>0$ such that $\Lambda_R\subset \res A_{\wedge,\Dom}$ and
\begin{equation*}
 \big\|\kappa_{|\lambda|^{1/m}}^{-1}(A_{\wedge,\Dom}-
  \lambda)^{-1}\big\|_{\L(x^{-m/2}L^2_b,\Dom_{\wedge,\max})}
  \leq C/|\lambda|,\quad \lambda\in \Lambda_R.
\end{equation*}
It can be shown that this estimate is equivalent to 
\eqref{NewFormWedgeCond} and \eqref{NecAndSuffWedgeCond}.

\begin{ex}
We consider again the model Laplacian $\Delta_{\wedge}$ from the previous 
section. Recall that $\bgres \Delta_{\wedge}=\C\backslash \overline{\R}_+$.  
For $\lambda\in \C\backslash\overline{\R}_+$, we have
\begin{equation*}
 \pi_{\wedge,\max}\K_{\wedge,\lambda}=
 \LinSpan\left\{-k_0 \log(-\lambda)+k_1\log x\right\}
 \;\text{ for some } k_0, k_1>0,
\end{equation*}  
where $\log$ means the principal branch of the logarithm. Moreover, 
by \eqref{FinitekappaD2}, 
\begin{equation*}
 \pi_{\wedge,\max}\kappa_{\!\varrho}^{-1}\Dom
 =\LinSpan\left\{\big(\tfrac{\zeta_0}{\zeta_1}-\log\varrho\big)\cdot 1
  + \log x\right\}.
\end{equation*}
The projection in \eqref{NecAndSuffWedgeCond} can be computed
explicitly. Namely, if $u=\alpha_0+\alpha_1\log x\in\Sing_{\wedge,\max}
=\LinSpan\{1,\log x\}$ and $\lambda=\varrho^{m}\lambda_0$, then
\begin{equation}\label{DeltaKmaxProjection}
\hat\pi_{\K_{\wedge,\lambda_0}, \kappa_\varrho^{-1}\Dom}u
 =\frac{-\alpha_0+\alpha_1(\frac{\zeta_0}{\zeta_1}-\log\varrho)}%
  {k_0\log(-\lambda_0)+k_1(\frac{\zeta_0}{\zeta_1}-\log\varrho)}
  \left(-k_0\log(-\lambda_0)+k_1\log x\right).
\end{equation}

Let $\Lambda$ be a closed sector in $\C\backslash \R_+$ containing the 
half-plane $\{\Re\lambda<0\}$.  Since the family of projections 
\eqref{DeltaKmaxProjection} is bounded as $\varrho\to\infty$, uniformly 
for $|\lambda_0|=1$ in $\Lambda$, regardless of the specific choice of 
$\alpha_0$, $\alpha_1$, Theorem~\ref{NecAndSuffWedge} implies that every 
closed extension $\Delta_{\wedge,\Dom}$, $\Dom\in\mathfrak G_\wedge$, of 
the model Laplacian admits $\Lambda$ as a sector of minimal growth.
\end{ex}

\subsection*{Equivalent geometric condition}
We identify $\mathfrak G_{\wedge}$ with the Grassmannian 
$\Gr_{d''}(\Sing_{\wedge,\max})$ where $d''=-\Ind(A_{\wedge,\min}-\lambda)$ 
for $\lambda\in\open\Lambda_\alpha\subset \bgres A_\wedge$. 
Let $d'=\dim\K_{\wedge,\lambda}$. The condition that in the Grassmannian 
$\Gr_{d''}(\Sing_{\wedge,\max})$, the curve 
\begin{equation*}
[R,\infty)\ni \varrho\mapsto \pi_{\wedge,\max}\kappa_\varrho^{-1}\Dom 
\end{equation*}
does not approach the set 
\[ \mathscr V_{\K_{\wedge,\lambda}}=\{D\in \Gr_{d''}(\Sing_{\wedge,\max}): 
   D\cap \pi_{\wedge,\max}\K_{\wedge,\lambda}\ne 0\} \]
as $\varrho\to\infty$, is sufficient for the validity of 
\eqref{NecAndSuffWedgeCond}. This is \cite[Theorem~8.28]{GKM1}. 
The following theorem states that the condition is also necessary.

For $D\in \Gr_{d''}(\Sing_{\wedge,\max})$ let
\begin{align*} 
\Omega^{-}(D)=\big\{D'\in \Gr_{d''}(\Sing_{\wedge,\max})
 :& \;\exists\, \{\varrho_k\}_{k=1}^\infty\subset \R_+ \text{ such that}\\ 
 &\; \varrho_k\to\infty \text{ and } 
 \kappa_{\varrho_k}^{-1}D\to D' \text{ as } k\to\infty\big\}. 
\end{align*}

\begin{thm}\label{GeometricNecAndSuff}
Let $\lambda_0\in\open\Lambda_\alpha$. The ray through $\lambda_0$ is a
ray of minimal growth for $A_{\wedge,\Dom}$ if and only if 
$\Omega^{-}(\pi_{\wedge,\max}\Dom) \cap \mathscr V_{\K_{\wedge,\lambda_0}} 
=\varnothing$.
\end{thm}
\begin{proof}
Let $\lambda_0\in\open\Lambda_\alpha$ and $\Dom\in\mathfrak G_\wedge$. 
For simplicity, we use the notation 
\begin{gather*} 
D=\pi_{\wedge,\max}\Dom, \quad 
  \mathscr V=\mathscr V_{\K_{\wedge,\lambda_0}}, \quad
  \K=\pi_{\wedge,\max}\K_{\wedge,\lambda_0} \;\text{ and }\; 
  \pi_{\K,D}=\hat\pi_{\K_{\wedge,\lambda_0},\Dom}.
\end{gather*}
Suppose $\Omega^{-}(D)\cap \mathscr V=\varnothing$. Since
$\Omega^{-}(D)$ and $\mathscr V$ are closed sets, there are a
neighborhood $\mathcal U$ of $\mathscr V$ and a constant $R>0$ 
such that if $\varrho>R$ then $\kappa_\varrho^{-1}D\not\in\mathscr V$.
Then Lemma~5.24 in \cite{GKM1} gives that
$\big\|\pi_{\K,\kappa_\varrho^{-1}D}\big\|$ is
uniformly bounded as $\varrho\to\infty$, and therefore, by
Theorem~\ref{NecAndSuffWedge} the ray through $\lambda_0$ is a ray of
minimal growth for $A_{\wedge,\Dom}$.

Assume now that there are $C$, $R>0$ such that $\Lambda_R\subset \res
A_{\wedge,\Dom}$ and the condition \eqref{NecAndSuffWedgeCond} is
satisfied. Suppose $\Omega^{-}(D)\cap \mathscr V \not=\varnothing$ 
and let $D_0\in \Omega^{-}(D)\cap \mathscr V$. Since $D_0\in
\mathscr V$, we have $D_0\cap\K\not=\{0\}$. On the other hand, $D_0\in 
\Omega^{-}(D)$ implies that there is a sequence 
$\{\varrho_k\}_{k=1}^\infty\subset\R_+$ such that $\varrho_k\to\infty$ and 
$D_k=\kappa_{\varrho_k}^{-1}D\to D_0$ as $k\to\infty$. Note that for 
$\varrho_k$ large we have $\varrho_k^m\lambda_0 \in\res A_{\wedge,\Dom}$, 
so $\lambda_0\in \res A_{\wedge,\kappa_{\varrho_k}^{-1}\Dom}$ and therefore,
$D_k\not\in \mathscr V$. 

Pick $v\in D_0\cap \K$ with $\|v\|=1$. Let $\pi_{D_k}$ denote the orthogonal 
projection on $D_k$. Since $D_k\to D_0$ as $k\to\infty$, we have $\pi_{D_k}
\to \pi_{D_0}$, so $v_k=\pi_{D_k}v\to \pi_{D_0}v=v$ as $k\to\infty$. Since 
$D_k\not\in \mathscr V$, $v_k-v\not=0$ and $\pi_{\K,D_k}v_k=0$. Hence
\[ \pi_{\K,D_k}\bigg(\frac{v-v_k}{\|v-v_k\|}\bigg) =
   \frac{v}{\|v-v_k\|} \to \infty \text{ as } k\to\infty, \]
since $\|v\|=1$ and $v_k\to v$ as $k\to\infty$. But this implies that
$\|\pi_{\K,D_k}\|\to\infty$ contradicting the boundedness of the norm
in \eqref{NecAndSuffWedgeCond}. Thus $\Omega^{-}(D)\cap 
\mathscr V =\varnothing$.
\end{proof}

\begin{ex}
Let $\Delta_{\wedge}$ be the model Laplacian and let 
$\Dom\in\mathfrak G_{\wedge}$. In this case, the limiting set 
$\Omega^{-}(\pi_{\wedge,\max}\Dom)$ consists of the one element of
$\mathbb{CP}^1$ corresponding to the Friedrichs extension of 
$\Delta_{\wedge}$, cf. \eqref{DomainConvergence}. From this new 
perspective, it is evident that every closed extension of 
$\Delta_{\wedge}$ must admit a sector of minimal growth.  
\end{ex}

\section{Rays of minimal growth}\label{sec-Necessity}

We continue to assume that $A\in x^{-m}\Diff_b^m(M;E)$ is $c$-elliptic.

Unlike the case of a differential operator with smooth coefficients on a 
closed manifold, that a ray $\Gamma$ is a ray of minimal growth for the 
principal symbol $\csym(A)$ of $A$ is not expected to imply that $\Gamma$ 
is a ray of minimal growth for $A$. In this context, it is useful to think 
of $A_\wedge$ as a symbol (the wedge symbol) associated with $A$, cf. 
Schulze \cite{Sz89}, so that it is natural to impose ray conditions on 
$A_\wedge$.  For this to work, however, we need to transfer the information 
about the given domain $\Dom$ of $A$ on $M$ to equivalent information for 
$A_\wedge$ on $Y^\wedge$, and vice versa.

\begin{thm}[Theorem~{4.12} in \cite{GKM1}] \label{thetaIsomorphism}
There is a natural isomorphism
\begin{equation*}
 \theta^{-1}:\Dom_{\wedge,\max}/\Dom_{\wedge,\min}\to \Dom_{\max}/\Dom_{\min}
\end{equation*}
given by a finite iterative procedure that involves the boundary spectrum
of $A$ and the decomposition  \eqref{bOperatorTaylor}. 
In particular, if $A$ has coefficients independent of $x$ near $Y$, then 
$\theta$ is the identity map.
\end{thm}

\begin{ex}
Let $M$ be a compact $2$-manifold with boundary $Y=S^1$. 
Let $A$ be a differential operator in $x^{-2}\Diff^2_b(M)$ that
over the interior of $M$ coincides with some Laplacian, and near 
$Y$, is of the form
\begin{equation*}
 A = x^{-2}\big((xD_x)^2 + q(x)\Delta_{Y}\big),  
\end{equation*}
where $\Delta_Y$ is the standard nonnegative Laplacian on $S^1$ and
$q$ is a smooth function. We assume $q$ to have the form
\begin{equation*}
  q(x)= \alpha^2 + \beta x + x^2\gamma(x),	
\end{equation*}
where $\alpha$, $\beta$ are constants such that $\frac12<\alpha<1$, 
$\beta\not=0$, and $\gamma(0)=1$. The associated model operator is
then given by 
\begin{equation*}
 A_\wedge = x^{-2}\big((xD_x)^2 + \alpha^2\Delta_{Y}\big),  
\end{equation*}
and $\spec_b(A)=\{\pm i\alpha k: k\in\mathbb N_0\}$. 
Since $\frac12<\alpha<1$, only the set $\{-i\alpha,0,i\alpha\}$ 
is relevant for the spaces $\Sing_{\max}$ and $\Sing_{\wedge,\max}$, cf. 
\eqref{SingularFunctions}. Here, similar to $\Sing_{\wedge,\max}$, the 
space $\Sing_{\max}$ consists of singular functions and is isomorphic to 
the quotient $\Dom_{\max}/\Dom_{\min}$.  
If $y$ denotes the angular variable on $S^1$, 
\begin{gather*}
\Sing_{\wedge,\max}=
\LinSpan\{1,\log x, e^{iy} x^{\alpha}, e^{-iy} x^{\alpha},
e^{iy} x^{-\alpha},e^{-iy} x^{-\alpha}\}, \\
\Sing_{\max}=
\LinSpan\Big\{1,\log x, e^{\pm iy} x^{\alpha}, 
e^{\pm iy} x^{-\alpha}\big(1-\tfrac{\beta}{2\alpha-1}x\big)\Big\}. 
\end{gather*}
In this case, $\theta:\Sing_{\max}\to \Sing_{\wedge,\max}$ acts as the 
identity on $\LinSpan\{1,\log x, e^{\pm iy} x^{\alpha}\}$, but
\begin{equation*}
 \theta\Big(e^{\pm iy}x^{-\alpha}
 \big(1-\tfrac{\beta}{2\alpha-1}x\big)\Big) = e^{\pm iy}x^{-\alpha}.
\end{equation*}
\end{ex}

\bigskip
The map $\theta$ induces an isomorphism
\begin{equation*}
  \Theta:\mathfrak{D} \to \mathfrak{D}_\wedge
\end{equation*}  
that we use to define $\Dom_\wedge=\Theta\Dom$ for any given
$\Dom\in \mathfrak{D}$. The operator $A_{\wedge,\Dom_\wedge}$ is the 
closed extension of $A_\wedge$ in $x^{-m/2}L_b^2(Y^\wedge;E)$ uniquely 
associated with $A_\Dom$.

As in \cite[Section 6]{GKM2}, and motivated by the importance of 
$\kappa_\varrho$ in studying the model operator $A_\wedge$, we introduce 
on $\Dom_{\max}(A)/\Dom_{\min}(A)$ the one-parameter group 
\begin{equation*}
 \tilde \kappa_\varrho = \theta^{-1} \kappa_\varrho \theta
 \;\text{ for } \varrho>0.
\end{equation*}

Similar to the situation on the model cone, the spectrum and resolvent of
the closed extensions of $A$ can be geometrically analyzed by considering 
the manifold $\mathfrak G$, cf. \eqref{GrasmannianA}, together with the 
flow generated by $\tilde\kappa_{\varrho}$.

An interesting consequence of Theorem~\ref{thetaIsomorphism} is the following.

\begin{prop}\label{IndexAwedgeD}
If $A-\lambda$ is $c$-elliptic with parameter $\lambda\not=0$, then
\[ \Ind(A_{\wedge,\Dom_\wedge}-\lambda)=\Ind A_\Dom. \]
\end{prop}
\begin{proof}
The existence of $\theta$ implies $\dim \Dom_\wedge/\Dom_{\wedge,\min}=
\dim \Dom/\Dom_{\min}$. Now, the proposition follows by combining this 
identity with the relative index formulas \eqref{RelIndexA} and 
\eqref{RelIndexAwedge}, together with the equation \eqref{IndexAwedgemin}.
\end{proof}

The following theorem describes the pseudodifferential structure of the 
resolvent of a cone operator $A$ and gives tangible conditions over a
given sector $\Lambda$ on the symbols $\csym(A)$ and $A_\wedge$ for $A$ to 
have $\Lambda$ as a sector of minimal growth.

\begin{thm}[Theorem~{6.9} in \cite{GKM2}]\label{MainPaper2}
Let $A\in x^{-m}\Diff_b^m(M;E)$ be such that $A-\lambda$ is $c$-elliptic 
with parameter $\lambda\in \Lambda$. If $\Lambda$ is a sector of minimal 
growth for $A_{\wedge,\Dom_\wedge}$, then it is a sector of minimal growth 
for $A_\Dom$. Moreover,
\[ (A_\Dom-\lambda)^{-1}=B(\lambda) + G_{\Dom}(\lambda), \]
where $B(\lambda)$ is a parametrix of $A_{\Dom_{\min}}-\lambda$ with
$B(\lambda)(A_{\Dom_{\min}}-\lambda)=1$ for $\lambda$ sufficiently large, 
and $G_{\Dom}(\lambda)$ is a pseudodifferential regularizing operator of 
finite rank.
\end{thm}

The following lemma gives further information about the behavior at large 
of the resolvent along a sector of minimal growth.

Given two cut-off functions $\omega_0$ and $\omega_1$, the
notation $\omega_1 \prec \omega_0$ will indicate that $\omega_0=1$ in 
a neighborhood of the support of $\omega_1$.

\begin{lem}\label{localizations}
Let $A \in x^{-m}\Diff_b^m(M;E)$ be $c$-elliptic and let $\Lambda$ be a 
sector of minimal growth for $A_{\Dom}$. For every pair of cut-off
functions $\omega_1\prec\omega_0$, supported near the boundary, we have
\begin{equation*}
(1-\omega_0)(A_{\Dom}-\lambda)^{-1}\omega_1
\in \S\bigl(\Lambda,\L(x^{-m/2}L_b^2,\Dom_{\max})\bigr),
\end{equation*}
where $\S$ stands for Schwartz (rapidly decreasing as $|\lambda|\to\infty$).
\end{lem}
\begin{proof}
Since $\Lambda$ is a sector of minimal growth for $A_{\Dom}$, the family
$A-\lambda$ must be $c$-elliptic with parameter $\lambda\in \Lambda$, and 
$A_{\wedge,\min}-\lambda$ must be injective for every
$\lambda\in\Lambda$, $\lambda\not=0$. 
A proof of this can be found in \cite[Theorem~4.1]{GKM2}.

As a consequence (cf. \cite[Section 5]{GKM2}), there is a parametrix 
$B(\lambda)$ such that $B(\lambda)(A_{\Dom_{\min}}-\lambda)=1$ for large 
$\lambda\in \Lambda$, and 
\begin{equation}\label{LocalizedParam}
(1-\omega_0) B(\lambda) \omega_1 \in 
\S\bigl(\Lambda,\L(x^{-m/2}L_b^2,\Dom_{\max})\bigr)
\end{equation}
for all cut-off functions $\omega_1\prec\omega_0$ supported near the
boundary. We now make use of the identity
\begin{equation*}
(A_{\Dom}-\lambda)^{-1} = B(\lambda) + 
(1-B(\lambda)(A-\lambda))(A_{\Dom}-\lambda)^{-1}. 
\end{equation*}
Multiplying by $(1-\omega_0)$ from the left and by $\omega_1$ from the 
right, \eqref{LocalizedParam} proves the assertion for the first term 
involving $B(\lambda)$. On the other hand, since $1-B(\lambda)(A-\lambda)$ 
vanishes on $\Dom_{\min}$ for large $\lambda$, we have for such $\lambda$,
\begin{align*}
(1-\omega_0)(1-B(\lambda)(A-\lambda)) 
 &= (1-\omega_0)(1-B(\lambda)(A-\lambda))\omega_2 \\
 &= -(1-\omega_0)B(\lambda)(A-\lambda)\omega_2 \\
 &= -(1-\omega_0)B(\lambda)\omega_1(A-\lambda)\omega_2
\end{align*}
whenever $\omega_2\prec\omega_1$.  Thus, by \eqref{LocalizedParam},
\begin{equation*}
(1-\omega_0)(1-B(\lambda)(A-\lambda)): \Dom_{\max} \to \Dom_{\max}
\end{equation*}
is rapidly decreasing as $|\lambda| \to \infty$. Finally, the assertion of 
the lemma can be completed using the fact that $(A_\Dom-\lambda)^{-1}
\omega_1: x^{-m/2}L_b^2 \to \Dom_{\max}$ is uniformly bounded. 
\end{proof}

\subsection*{Necessity of the conditions}

The converse of Theorem~\ref{MainPaper2} involves proving that the 
minimal growth of the resolvent $(A_\Dom-\lambda)^{-1}$ over a sector 
$\Lambda$ implies a corresponding behavior for the inverse of $\csym(A)-
\lambda$ and for the resolvent $(A_{\wedge,\Dom_{\wedge}}-\lambda)^{-1}$. 

While in \cite[Theorem~4.1]{GKM2} we established the necessity of the
condition on $\csym(A)$, we did not address the question whether
$\Lambda$ must necessarily be a sector of minimal growth for 
$A_{\wedge,\Dom_{\wedge}}$.  In the next theorem we prove that this is 
indeed the case when $A$ has coefficients independent of $x$ near 
$Y=\partial M$.

\begin{thm}\label{NecessityConstCoeff}
Let $A\in x^{-m}\Diff_b^m(M;E)$ be $c$-elliptic with coefficients 
independent of $x$ near $Y$. If $\Lambda$ is a a sector of minimal growth 
for $A_\Dom$, then $A-\lambda$ is $c$-elliptic with parameter 
$\lambda\in\Lambda$, and $\Lambda$ is a sector of minimal growth for 
$A_{\wedge,\Dom_\wedge}$.
\end{thm}

\begin{proof}
As stated in the proof of Lemma~\ref{localizations}, the assumption on
the resolvent of $A_{\Dom}$ implies that $A-\lambda$  is $c$-elliptic with 
parameter $\lambda\in \Lambda$ and that $A_{\wedge,\min}-\lambda$ is 
injective for every $\lambda\not=0$. Thus we only need to prove the
statement about $A_{\wedge,\Dom_\wedge}$.

By Proposition~\ref{IndexAwedgeD}, and since $\Ind A_\Dom=0$,
we have $\Ind(A_{\wedge,\Dom_\wedge}-\lambda)=0$ for $\lambda\not=0$. 
For this reason, in order to show that $\Lambda$ is a sector of minimal 
growth for $A_{\wedge,\Dom_\wedge}$, it suffices to find (for large 
$\lambda\in\Lambda$) a right-inverse of $A_{\wedge,\Dom_\wedge}-\lambda$ 
that is uniformly bounded in $\L\big(x^{-m/2}L_b^2,\Dom_\wedge\big)$ as 
$|\lambda|\to\infty$.

Since $A$ is assumed to have coefficients independent of $x$ near the
boundary, there is a cut-off function $\omega_0$ such that
\begin{equation*}
 A\omega_0=A_\wedge \omega_0 \quad\text{and}\quad  
 \omega_0\Dom=\omega_0\Dom_\wedge. 
\end{equation*}
Let $\omega_1$, $\omega_2$ be cut-off functions with $\omega_2\prec
\omega_1\prec \omega_0$. Then the operator 
\begin{equation*}
 B(\lambda)=\omega_1 (A_\Dom-\lambda)^{-1} \omega_2
\end{equation*}
can be regarded as an operator on $M$ with values in $\Dom$ or as an
operator on $Y^\wedge$ with values in $\Dom_\wedge$.  Depending on the
context we will write $B(\lambda)$ as
\begin{equation*}
B_\Dom(\lambda):x^{-m/2}L^2_b(M;E)\to\Dom \quad\text{or}\quad 
B_{\Dom_\wedge}(\lambda):x^{-m/2}L^2_b(Y^\wedge;E)\to\Dom_{\wedge}.
\end{equation*}
On $M$ we consider
\begin{align*}
(A_\Dom-\lambda)B_\Dom(\lambda)
&=\omega_0(A_\Dom-\lambda)\omega_1 (A_\Dom-\lambda)^{-1} \omega_2 \\
&=\omega_2 -\omega_0(A_\Dom-\lambda)(1-\omega_1)
  (A_\Dom-\lambda)^{-1}\omega_2 \\
&=\omega_2 + R(\lambda)
\end{align*}
with $R(\lambda)=-\omega_0(A_\Dom-\lambda)(1-\omega_1)(A_\Dom-\lambda)^{-1}
\omega_2$. By Lemma~\ref{localizations},  $R(\lambda)$ is rapidly
decreasing in the norm as $|\lambda|\to \infty$.

Because of the presence and nature of the cut-off functions $\omega_0$ and 
$\omega_2$, $R(\lambda)$ can also be regarded as an operator on $Y^\wedge$, 
say $R_\wedge(\lambda)\in \S\big(\Lambda,\L(x^{-m/2}L^2_b)\big)$. 
Now, using that $(A_\Dom-\lambda)\omega_1= (A_{\wedge,\Dom_\wedge}-\lambda)
\omega_1$, we get on $Y^\wedge$ the identity
\begin{equation}\label{ParamNearBoundary}
 (A_{\wedge,\Dom_{\wedge}}-\lambda)B_{\Dom_\wedge}(\lambda)=
 \omega_2 + R_\wedge(\lambda).
\end{equation}
Furthermore, we have
\begin{equation*}
\|B_{\Dom_\wedge}(\lambda)\|_{\L(x^{-m/2}L^2_b, 
\Dom_{\wedge,\max})}= O(1) \;\text{ as } |\lambda|\to\infty, 
\end{equation*}
since, by assumption, $\|B_{\Dom}(\lambda)\|_{\L(x^{-m/2}L^2_b,\Dom_{\max})}$ 
has the same asymptotic behavior.

On the other hand, as $A-\lambda$ is $c$-elliptic with parameter, by
\cite[Theorem~5.24]{GKM2} there is a family of pseudodifferential
operators $B_{2,\wedge}(\lambda):x^{-m/2}L^2_b\to \Dom_{\wedge,\min}$ 
(uniformly bounded in $\lambda$) such that $(A_\wedge-\lambda)
B_{2,\wedge}(\lambda)-1$ is regularizing, and for $\omega_3\prec \omega_2$,
the families $\omega_3 B_{2,\wedge}(\lambda) (1-\omega_2)$ and
$\big[(A_\wedge-\lambda)B_{2,\wedge}(\lambda)-1\big](1-\omega_2)$ are
rapidly decreasing in the norm as $|\lambda| \to \infty$. Thus, as
$A_{\wedge,\Dom_\wedge}(1-\omega_3)=A_{\wedge}(1-\omega_3)$, 
\begin{equation}\label{ParamAwayBoundary}
(A_{\wedge,\Dom_\wedge}-\lambda)(1-\omega_3)B_{2,\wedge}(\lambda)(1-\omega_2)
= (1-\omega_2) + S_\wedge(\lambda)
\end{equation}
with $S_\wedge(\lambda)\in \S\big(\Lambda,\L(x^{-m/2}L^2_b)\big)$.
Finally, the operator family
\begin{equation*}
Q_\wedge(\lambda)
=B_{\Dom_\wedge}(\lambda)+(1-\omega_3)B_{2,\wedge}(\lambda)(1-\omega_2)
:x^{-m/2}L^2_b\to \Dom_{\wedge,\max}
\end{equation*}
is bounded in the norm as $|\lambda|\to\infty$ and by 
\eqref{ParamNearBoundary} and \eqref{ParamAwayBoundary} we have
\begin{equation*}
(A_{\wedge,\Dom_\wedge}-\lambda)Q_\wedge(\lambda)-1 \in
\S\big(\Lambda,\L(x^{-m/2}L^2_b)\big).
\end{equation*}
By a Neumann series argument, it follows that
$A_{\wedge,\Dom_\wedge}-\lambda:\Dom_\wedge\to x^{-m/2}L^2_b$ has a
uniformly bounded right-inverse for large $\lambda\in\Lambda$.
\end{proof}



\subsection*{Acknowledgment}
The new results contained herein reflect part of work carried out by the 
three authors at the Mathematisches Forschungsinstitut Oberwolfach under 
their ``Research in Pairs'' program. They gratefully acknowledge the 
Institute's support and hospitality. 


\begin{thebibliography}{99}

\bibitem{AV63}
M.~Agranovich and M. Vishik, \emph{Elliptic problems with a
parameter and parabolic problems of general type}, Russ. Math. Surveys
\textbf{19} (1963), 53--159.

\bibitem{Agmon}
S.~Agmon, \emph{On the eigenfunctions and on the eigenvalues of general 
elliptic boundary value problems}, Comm. Pure Appl. Math. \textbf{15}
(1962), 119--147.

\bibitem{BrSe91}
J.~Br\"uning and R.~Seeley, \emph{The expansion of the resolvent
near a singular stratum of conical type},
J. Funct. Anal. \textbf{95} (1991), no. 2, 255--290.

\bibitem{Gil}
J.~Gil, \emph{Full asymptotic expansion of the heat trace for
non-self-adjoint elliptic cone
operators}, Math. Nachr. \textbf{250} (2003), 25--57.

\bibitem{GKM1}
J.~Gil, T.~Krainer, and G.~Mendoza, \emph{Geometry and spectra of 
closed extensions of elliptic cone operators}, preprint math.AP/0410178 
at arXiv.org, to appear in Canadian Journal of Mathematics.

\bibitem{GKM2}
J.~Gil, T.~Krainer, and G.~Mendoza, \emph{Resolvents of
elliptic cone operators}, preprint math.AP/0410176 at arXiv.org.

\bibitem{GiMe01}
J.~Gil and G.~Mendoza, \emph{Adjoints of elliptic cone operators},
Amer. J. Math. \textbf{125} (2003) 2, 357--408.

\bibitem{GrubbBuch} 
G.~Grubb, \emph{Functional calculus of pseudodifferential boundary problems},
2nd ed., Progress in Mathematics, vol.~65. Birkh{\"a}user, Basel, 1996.

\bibitem{K1} 
V.~Kondrat'ev, \emph{Boundary problems for
elliptic equations in domains with conical or angular points},
Trans. Mosc. Math. Soc. \textbf{16} (1967), 227--313.

\bibitem{Krainer}
T.~Krainer, \emph{Resolvents of elliptic boundary problems on conic
manifolds}, preprint math.AP/0503021 at arXiv.org.

\bibitem{Le97} 
M.~Lesch, \emph{Operators of {F}uchs type, conical
singularities, and asymptotic methods}, Teubner-Texte zur Math.
vol 136, B.G. Teubner, Stuttgart, Leipzig, 1997.

\bibitem{LoRes01}
P.~Loya, \emph{On the resolvent of differential operators on conic
manifolds}, Comm. Anal. Geom. \textbf{10} (2002), no.~5, 877--934.

\bibitem{RBM1}
R.~Melrose, \emph{Transformation of boundary value
problems}, Acta Math. \textbf{147} (1981), no. 3-4, 149--236.

\bibitem{SchrSe03}
E. Schrohe and J. Seiler,
\emph{The resolvent of closed extensions of cone differential
operators}, Canad. J. Math. \textbf{57} (2005), no. 4, 771--811.

\bibitem{Sz89}
B.-W.~Schulze, \emph{Pseudo-differential operators
on manifolds with edges}, Proc. Symp. `Partial Differential
Equations', Holzhau 1988 (Leipzig), Teubner-Texte zur Math., vol.
112, Teubner, 1989, pp.~259--288.

\bibitem{Seeley}
R.~Seeley, \emph{Complex powers of an elliptic operator}, Singular
Integrals, AMS Proc. Symp. Pure Math. X, 1966, Amer. Math. Soc.,
Providence, 1967, pp.~288--307.

\end{thebibliography}
\end{document}